\documentclass{amsart}
\usepackage[final]{epsfig}
\usepackage{graphics}
\usepackage{float}
\usepackage{amsmath}
\usepackage{amsfonts}
\usepackage{latexsym}
\usepackage{amssymb}
\usepackage{graphicx}
\usepackage{url}
\usepackage{epstopdf}
\usepackage{verbatim}
\usepackage[margin=1in]{geometry}
\usepackage{amsthm}
\usepackage{tikz}
\usepackage{caption}

\makeatletter
\newtheorem*{rep@theorem}{\rep@title}

\newcommand{\newreptheorem}[2]{%
\newenvironment{rep#1}[1]{%
 \def\rep@title{#2 \ref{##1}}%
 \begin{rep@theorem}}%
 {\end{rep@theorem}}}
\makeatother

\newtheorem{lemma}{Lemma}[section]

\newtheorem{example}[lemma]{Example}
\newtheorem{theorem}[lemma]{Theorem}

\newtheorem{corollary}[lemma]{Corollary}

\newtheorem*{theorem*}{Theorem}

\newreptheorem{theorem}{Theorem}

\newcommand{\proofend}{$\Box$\bigskip}

\newcommand{\N}{{\mathbb N}}

\newcommand{\Z}{{\mathbb Z}}

\newcommand{\lfl}{\lfloor}
\newcommand{\rfl}{\rfloor}
\newcommand{\lrep}{\left[}
\newcommand{\rrep}{\right]}
\def\proof{\par{\it Proof}. \ignorespaces}

\DeclareMathOperator{\Aut}{Aut}

\begin{document}

\title{Equality of Dedekind sums mod $\mathbb{Z},2\mathbb{Z}$ and $4\mathbb{Z}$}

\author{E. Tsukerman}
\date{\today}
\address{Department of Mathematics, University of California,
Berkeley, CA 94720-3840}
\email{e.tsukerman@berkeley.edu}
\subjclass[2000]{}

\begin{abstract}

In [Girstmair, A criterion for the equality of Dedekind sums mod $\mathbb{Z}$, Internat. J. Number Theory 10: (2014) 565--568], it was shown that the necessary condition $b \mid (a_1 a_2-1)(a_1-a_2)$ for equality of two dedekind sums $s(a_1,b)$ and $s(a_2,b)$ given in [Jabuka, Robins and Wang, When are two Dedekind sums equal? Internat. J. Number Theory 7: (2011) 2197--2202] is equivalent to $12s(a_1,b)-12s(a_2,b) \in \mathbb{Z}$. In this note, we give a new proof of this result and then find two additional necessary and sufficient conditions for $12s(a_1,b)-12s(a_2,b) \in  2\mathbb{Z}, 4\mathbb{Z}$. These give new necessary conditions on equality of Dedekind sums.
\end{abstract}

\maketitle

In \cite{MR2873148}, Jabuka et al. raise the question of when two Dedekind sums $s(a_1,b)$ and $s(a_2,b)$ are equal. In the same paper, they prove the necessary condition $b \mid (a_1 a_2-1)(a_1-a_2)$ for equality of two dedekind sums $s(a_1,b)$ and $s(a_2,b)$. Girstmair \cite{MR3189994} shows that this condition is equivalent to $12s(a_1,b)-12s(a_2,b) \in \mathbb{Z}$. In this note, we give a new proof of this result and then find two additional necessary and sufficient conditions for $12s(a_1,b)-12s(a_2,b) \in  2\mathbb{Z}, 4\mathbb{Z}$. These give new necessary conditions on equality of Dedekind sums.

\section{Preliminaries}

Let  $ \pi_{(a,b)} \in \Aut(\Z / b \Z), \pi_{(a,b)}:x \mapsto ax$. Let $\lrep x \rrep_b=x-b \lfl \frac{x}{b} \rfl$ be the function taking $x \in \Z/ b \Z$ to its smallest nonnegative representative. We view $\pi_{(a,b)}$ as a permutation of $\{0,1,\ldots,b-1\}$ given by
\[
\pi_{(a,b)}=\left(\begin{array}{ccccccc}
0 & 1 & \cdots & b-1 \\
\lrep \pi(0) \rrep & \lrep \pi(1) \rrep & \cdots & \lrep \pi(b-1) \rrep
\end{array}\right)=\left(\begin{array}{ccccccc}
0 & 1 & \cdots & b-1 \\
0& \lrep a \rrep_b & \cdots & \lrep (b-1)a \rrep_b
\end{array}\right).
\]
The precedent for doing so is already present in the work of Zolotarev, in which he relates the sign of $\pi_{(a,b)}$ to the Jacobi symbol and obtains a proof of the law of quadratic reciprocity (see, e.g., \cite[pg. 38]{MR0357299}). Let $I(a,b)$ denote the number of inversions  of $\pi_{(a,b)}$.

\begin{theorem}  (Zolotarev) For odd $b$ and $(a,b)=1$,
\[
(-1)^{I(a,b)}=\left( \frac{a}{b} \right).
\]
\end{theorem}

The following result shows that the inversions of $\pi_{(a,b)}$ and Dedekind sums are closely related.

\begin{theorem} \label{invded}(Meyer, \cite{Meyer1957})  The number of inversions $I(a,b)$ of $\pi_{(a,b)}$ is equal to
\[
I(a,b)=-3b s(a,b)+\frac{1}{4}(b-1)(b-2),
\]
where $s(a,b)$ is the Dedekind sum.
\end{theorem}

From the reciprocity law of dedekind sums, one obtains a reciprocity law for inversions. 

\begin{theorem} \label{sal}(Sali\'{e})  For all coprime $a,b \in \N$
\[
4aI(a,b)+4bI(b,a)=(a-1)(b-1)(a+b-1).
\]
\end{theorem}

\section{Necessary Condition on Equality of Dedekind Sums}

The equivalence of 1(b) and 1(c) of Theorem \ref{necCond} is proved in \cite{MR3189994} via the Barkan-Hickerson-Knuth formula.

\begin{theorem} \label{necCond} Let $a_1,a_2 \in \N$ be relatively prime to $b \in \N$. 
\begin{enumerate}\label{1st}
\item The following are equivalent:
\begin{enumerate}
\item $4I(a_1,b) \equiv 4I(a_2,b) \pmod{b}$
\item $12s(a_1,b)-12s(a_2,b) \in \Z$
\item $b \mid (a_1-a_2)(a_1a_2-1)$.
\end{enumerate}
\item The following are equivalent:
\begin{enumerate}
\item $2I(a_1,b) \equiv 2I(a_2,b) \pmod{b}$
\item $6s(a_1,b)-6s(a_2,b) \in \Z$
\item $2b \mid (a_1-a_2)(b-1)(b+a_1 a_2-1)$.
\end{enumerate}
\item The following are equivalent:
\begin{enumerate}
\item $I(a_1,b) \equiv I(a_2,b) \pmod{b}$
\item $3s(a_1,b)-3s(a_2,b) \in \Z$
\item $4b \mid (a_1-a_2)(b-1)(b+a_1 a_2-1)$.
\end{enumerate}
\end{enumerate}
\end{theorem}

\proof
From Theorem \ref{invded}, it is immediate that
\begin{align}
4I(a_1,b)-4I(a_2,b)=-b(12s(a_1,b)-12s(a_2,b)).
\end{align}
This shows the equivalence of all $(a)$'s with corresponding $(b)$'s.

Suppose that $m \in \{b,2b,4b\}$. Reducing both sides of the equation in Theorem \ref{sal} mod $m$ yields
\[
4 a_1 I(a_1,b) \equiv (a_1-1)(b-1)(a_1+b-1) \pmod{m}.
\]
\[
4 a_2 I(a_2,b) \equiv (a_2-1)(b-1)(a_2+b-1) \pmod{m}.
\]
Multiplying the first congruence by $a_2$, the second by $a_1$, and subtracting, we obtain
\begin{align}\label{modm}
4a_1 a_2 (I(a_1,b)-I(a_2,b)) \equiv (a_1-a_2)(b-1)(b+a_1 a_2-1) \pmod{m}.
\end{align}
The equivalence of $(a)$'s and $(c)$'s follows from (\ref{modm}) and $\gcd(a_1,b)=\gcd(a_2,b)=1$.
\proofend

\begin{example}
We consider the example with $a_1=37, a_2=33$ and $b=40$ given in \cite{MR2873148}. Since 
\[
b=40 \mid (a_1-a_2)(a_1 a_2-1)=4880,
\]
we have $12s(a_1,b)-12s(a_2,b) \in \Z$. We also see that
\[
2b=80 \mid 4880,
\]
so that $12s(a_1,b)-12s(a_2,b) \in 2\Z$. However,
\[
4b=120 \nmid 4880,
\]
implying that $12s(a_1,b)-12s(a_2,b) \not \in 4\Z$. Indeed, $s(37,40)=-\frac{13}{16}$ and $s(33,40)=-\frac{5}{16}$. Thus,
\[
12s(33,40)-12s(37,40)=6.
\]
\end{example}

\begin{example}
Let $a_1=1, a_2=6$ and $b=25$. Then
\[
12s(1,25)-12s(6,25)=24.
\]
This shows that $a_1=1$ and $a_2=6$ satisfy the necessary conditions in \ref{necCond} but do not yield an equality of Dedekind sums.
\end{example}

\begin{corollary}\label{e}
Let $b=2^e p$ for $e \in \{0,1,2,3\}$ and $p$ an odd prime. Then $s(a_1,b)=s(a_2,b)$ if and only if $a_1 \equiv a_2 \pmod{b}$ or $a_1 \equiv a_2^{-1} \pmod{b}$.
\end{corollary}

\proof
If $e=0$, then $b$ is prime and the result clearly holds. Assume then that $e \geq 1$. We will employ the fact that since $a_1$ is odd and $e \leq 3$, $a_1^2 \equiv 1 \pmod{2^e}$. Assume that $s(a_1,b)=s(a_2,b)$, so that 
\[
4b \mid (a_1-a_2)(b-1)(b-1+a_1a_2).
\]
There are two cases to consider: when $p \mid (a_1-a_2)$ and when $p \mid (a_1 a_2-1)$. Consider the first case.
Write $a_1=a_2+pk$. Let $k=2^m k'$ with $k'$ odd. We also $a_2^2-1=2^e j$ for $j$ odd. Then
\[
2^{e+2} \mid 2^m k'(2^ep+2^ej+2^m a_2pk).
\]
Consequently, $e+2 \leq m+\min\{e,m\}$. The assumption that $e \leq 3$ implies that $m \geq e$.
For the second case, we write $a_2=\frac{1+pk}{a_1}$. Then
\[
2^{e+2} \mid (a_1-\frac{1+pk}{a_1})(2^ep+pk) \implies 2^{e+2} \mid (a_1^2-1-pk)(2^e p+pk).
\]
Set $k=2^m k'$ and $a_1^2=1+2^ej$. Then
\[
e+2 \leq 2 \min\{m,e\},
\]
and the result follows.
\proofend

\begin{example}
Corollary \ref{e} does not hold for $e>3$. Indeed, take $b=3 \cdot 2^4=48$, $a_1=1$ and $a_2=25$. Then
\[
s(1,48)=\frac{1081}{288}, \quad s(25,48)=\frac{217}{288},
\]
so that $12s(1,48)-12s(25,48)=36$.
\end{example}

\bigskip
{\bf Acknowledgments}.  
This material is based upon work supported by the National Science Foundation Graduate Research Fellowship under Grant No. DGE 1106400. Any opinion, findings, and conclusions or recommendations expressed in this material are those of the authors(s) and do not necessarily reflect the views of the National Science Foundation.

\bibliographystyle{alpha}
\bibliography{bibliography}

\end{document}